\documentclass[12pt]{article}
\usepackage{amsfonts, amsmath, amssymb,amsthm,bbm,comment,fancyhdr,float,makeidx,mathrsfs,verbatim}

\normalsize
\newcommand{\al}{\alpha}

\newcommand{\seq}{\subseteq}

\newcommand{\Lra}{\Longrightarrow}

\newcommand{\ol}{\overline}

\newcommand{\vs}{\vspace*}

\newcommand{\nin}{\noindent}

\newtheorem{mthm}{Theorem}[section]
\newtheorem{mylem}[mthm]{Lemma}
\newtheorem{myprn}[mthm]{Proposition}
\newtheorem{mycor}[mthm]{Corollary}
\newtheorem{mydef}[mthm]{Definition}
\newtheorem{myrem}[mthm]{Remark}
\newtheorem{mycon}[mthm]{Construction}
\newtheorem{myeg} [mthm]{Example}
\newtheorem{myque} [mthm]{Question}
\newenvironment{thm}{\begin{mthm}}{\end{mthm}}
\newenvironment{lem}{\begin{mylem}}{\end{mylem}}
\newenvironment{cor}{\begin{mycor}}{\end{mycor}}
\newenvironment{prop}{\begin{myprn}}{\end{myprn}}
\newenvironment{mdef}{\begin{mydef}\rm}{\end{mydef}}
\newenvironment{rem}{\begin{myrem}\rm}{\end{myrem}}

\newenvironment{ex}{\begin{myeg}\rm}{\end{myeg}}
\newenvironment{prof}{\noindent $Proof.$ \rm}{\hfill $\Box$}

\def\Lra{\Longrightarrow}
\def \nin {\noindent}

\def \Lemma #1 {\vs{3mm}\nin {\bf Lemma #1} \it}
\def \Prop #1 {\vs{3mm}\nin {\bf Proposition #1} \it}
\def \Th #1 {\vs{3mm}\nin {\bf Theorem #1} \it}
\def \Cor #1 {\vs{3mm}\nin {\bf Corollary #1} \it}
\def \Ex #1 {\vs{3mm}\nin {\bf Example #1} \it}

\def \part #1 {\hfil\break\hglue 12pt {\rm (#1)~}}

\def\fs{\footnotesize}

\title{
\bf\LARGE  Monomial ideals under ideal operations\thanks{This research is supported by the National Natural
Science Foundation of China (Grant No. 11271250). } }
\author{Jin Guo\thanks{ guojinecho@163.com} \,\,and {Tongsuo Wu\thanks{Corresponding author. tswu@sjtu.edu.cn}}\\
 {\small Department of Mathematics, Shanghai Jiaotong University}
}

\date{}
\begin{document}
\baselineskip=16pt \maketitle

\begin{center}
\begin{minipage}{12cm}

 \vs{3mm}\nin{\small\bf Abstract.} {\fs  In this paper, we show for a monomial ideal $I$ of $K[x_1,x_2,\ldots,x_n]$ that the integral closure $\ol{I}$ is a monomial ideal of Borel type  (Borel-fixed, strongly stable, lexsegment, or universal lexsegment respectively), if $I$ has the same property. We also show that the $k^{th}$ symbolic power $I^{(k)}$ of $I$ preserves the properties of Borel type, Borel-fixed and strongly stable, and $I^{(k)}$ is lexsegment if $I$ is stably lexsegment.  For a monomial ideal $I$ and a monomial prime ideal $P$, a new ideal $J(I, P)$ is studied, which also gives a clear description of the primary decomposition of $I^{(k)}$.  Then a new simplicial complex $_J\bigtriangleup$ of a monomial ideal $J$ is defined, and it is shown that $I_{_J\bigtriangleup^{\vee}} = \sqrt{J}$. Finally, we  show under an additional weak assumption that a monomial ideal is universal lexsegment if and only if its polarization is a squarefree strongly stable ideal.  }

\vs{3mm}\nin {\small Key Words:} {\small  Borel type monomial ideal; $k^{th}$ symbolic power; integral closure; polarization; universal lexsegment  monomial ideal}

\end{minipage}
\end{center}

\vs{4mm} \section{Introduction }

\vs{3mm}\nin   Throughout the  paper, $K$ is an infinite field and let $S=K[x_1,\,x_2,\,\ldots,\,x_n]$ be the polynomial ring with $n$ indeterminants over $K$. If an ideal $I$ is generated by $u_1,\ldots,u_s$, then we denote it by $I=\langle u_1,\ldots,u_s\rangle$. For a monomial ideal $I$ of $S$, recall that $I$ is called {\it strongly stable} if for any monomial $u$ in $I$ and any $i<j\le n$, $x_j\mid u$ implies $x_i(u/x_j)\in I$. Recall that $I$ is called {\it Borel-fixed}, if $\al(u)\in I$ holds for any invertible upper  $n\times n$ matrix $\al$ over $K$.  Recall that  $I$  is called {\it of Borel type} if
$$I:x_i^{\infty}=I: \langle x_1,\,x_2,\,\ldots,\,x_i \rangle^\infty\quad \quad (*)$$
holds for every  $i=1,\ldots,n$. It is known that each strongly stable monomial ideal is Borel-fixed, and the converse holds under the additional assumption $char(K)=0$.  Bayer and Stillman in \cite{BS} noted that Borel-fixed ideals satisfy condition $(*)$. Herzog et al. in \cite{HPV} gave the definition of  a  Borel type monomial ideal, and they proved among other things that  a  Borel type monomial ideal is sequentially Cohen-Macaulay, see also \cite{Popescu}. Furthermore, there are other two classes of strongly stable monomial ideals, namely, monomial ideals which are lexsegment or universal lexsegment,  see \cite{AHH} or \cite{HH}. We have the following relations for conditions on a monomial ideal:

\vs{3mm}

\noindent universal lexsegment$\Lra$lexsegment$\Lra$strongly stable$\Lra$ Borel-fixed$\Lra$ of Borel type.

\vs{2mm} The following is the fundamental characterization of Borel type monomial ideals:

\begin{prop}\label{BT} (\cite[Proposition 4.2.9]{HH}) For a monomial ideal $I$ of $S$, the following conditions are equivalent:

$(1)$ $I$ is of Borel type.

$(2)$ For each monomial $u\in I$ and all positive  integers $i,j,s$ with $i<j\le n$ such that $x_j^s\mid u$, there exists an integer $t\ge 0$ such that $x_i^t(u/x_j^s)\in I$.

$(3)$ Each associated prime ideal $P$ of $I$ has the form $\langle x_1,x_2,\ldots, x_r \rangle$ for some $r\le n$.
\end{prop}



In \cite[Proposition 1]{MC}, Mircea Cimpoeas observed that the afore mentioned property is preserved under several operations, such as sum, intersection, product, colon. For  a monomial ideal $I$ of Borel type, note that $I:\mathfrak{m}^\infty=I:\mathfrak{m}^r$ holds for $r>>0$, thus the saturation $I:\mathfrak{m}^\infty$ is a monomial ideal of Borel type. The root ideal $\sqrt{I}$ is a prime ideal of the form $\langle x_1,x_2,\ldots,x_r \rangle$, and is thus universal lexsegment.

Some parts of the following proposition are well known, the others are direct to check, so we omit the verification.

\begin{prop}\label{operation} Let $I,J, L$ be monomial ideals of $S$.

(1) If further $I,J$ are
of Borel type (strongly stable,  respectively), then each of the following is a monomial ideal of Borel type ( strongly stable, respectively):
$$ I\cap J ,\, I+J,\,\,I:L,\,\, IJ. $$
In particular, the saturation $I:\mathfrak{m}^\infty $ of $I$ is of Borel type (strongly stable, respectively) if $I$ has the same property.

(2) If further $I,J$ are Borel-fixed ideals, then each of $ I\cap J ,\, I+J,\,\, I:J,\,\, IJ$ is again Borel-fixed. In particular, the saturation $I:\mathfrak{m}^\infty $ of $I$ is Borel-fixed.

(3) If further $I,J$ are lexsegment (universal lexsegment, respectively) ideals, then each of $ I\cap J ,\, I+J,\,\, I:L$ is again lexsegment (universal lexsegment, respectively).
\end{prop}

\vs{3mm} Let $I$ be a Borel-fixed monomial ideal, and $L$  a monomial ideal which need not to be Borel-fixed. The following example shows that the colon $I : L$  may be not Borel-fixed.

\begin{ex}\label{colon not Borel-fixed} Let $K$ be a field with $char(K)=2$, and let $S=K[x_1, \ldots, x_n]$. If $I= \langle x_1^3, x_1x_2^2 \rangle$. It is direct to check that $I$ is Borel-fixed. Set $L = \langle x_2 \rangle$. It is easy to see that $I : L = \langle x_1^3, x_1x_2 \rangle$, which is not Borel-fixed.
\end{ex}

\vs{3mm} The following example shows that $IJ$ may be not lexsegment, even though $I, J$ are lexsegment.

\begin{ex}\label{product not lexsegment} Let $S=K[x_1, x_2, x_3]$, and let $I= \langle x_1^3, x_1^2x_2, x_1^2x_3, x_1x_2^2, x_1x_2x_3 \rangle$. It is easy to see that $I$ is lexsegment, and $u=x_1^2x_2^2x_3^2 \in I^2$. Note that $v= x_1^3x_3^3 \not\in I^2$ and $v >_{lex}u$, so $I^2$ is not lexsegment.
\end{ex}


\vs{3mm}As an application of Proposition \ref{operation}, we now give an alternative proof to the following:

\begin{cor}\label{regular} (\cite [Proposition 4.3.3]{HH})
Let $I\seq S$ be a monomial ideal of Borel type. Then $x_n,\ldots,x_1$ is an almost regular sequence on $S/I$.
\end{cor}

\begin{prof}  In the proof of \cite [Lemma 4.3.1]{HH}, let $M=S/I$. Then the corresponding $N$ (i.e.,  $0:_M\mathfrak{m}^\infty$) is  identical with $(I:\mathfrak{m}^\infty)/I$. Note that
$M/N\cong S/(I: \mathfrak{m}^\infty )$ holds. If $M=N$, then each element of $S_1$ is almost regular on $M$ since $M$ has finite length. Now assume $M\neq N$. Since $I: \mathfrak{m}^\infty $ is monomial of Borel type and $\mathfrak{m}\not\in Ass(M/N)$, as is shown in the proof of the Lemma 4.3.1,  it follows by \cite [Proposition 4.2.9(d)]{HH} that $x_n\not\in \cup Ass(M/N)$, i.e., $x_n$ is in the constructed open set $U$ and thus is almost regular on $S/I$.
The result then follows by mathematical induction.
\end{prof}




\section{Integral Closure $\ol{I}$}

Let  $I$ be any ideal of a commutative ring $R$. Recall from \cite{Swanson} that the integral closure $\ol{I}$ of an ideal $I$  consists of all elements of $R$ which are integral
over $I$. Note that $\ol{I}$ is an ideal of $R$. For a monomial ideal $I$ of $S$, $\ol{I}$ is generated  by all monomials $u$ such that $u^k\in I^k$ holds for some $k>0$. Thus the exponent set of all monomials in $\ol{I}$ is identical with the integer lattice points in the convex hull of the exponent set of all monomials in $I$, see \cite[Proposition 1.4.6]{Swanson}. In this section, we will show that $\ol{I}$ is a monomial ideal of Borel type (strongly stable, Borel-fixed, lexsegment, or universal lexsegment  respectively), whenever $I$ has the same property.




\begin{thm}\label{Closure} Let $I$ be a monomial ideal of $S$ and let $\ol{I}$ be its integral closure. If $I$ is of Borel type (strongly stable, Borel-fixed, lexsegment, or universal lexsegment respectively), then $\ol{I}$ is also monomial of Borel type (strongly stable, Borel-fixed, lexsegment, or universal lexsegment  respectively).

\end{thm}

\begin{prof}
(1) Assume that $I$ is monomial of Borel type. Then $\ol{I}$ is also monomial by \cite[Proposition 1.4.2]{Swanson} (see also \cite[Theorem 1.4.2]{HH}). In order to prove that $\ol{I}$ is of Borel type, we need only to verify that each associated prime ideal of $\ol{I}$ has the form $\langle x_1,x_2,\ldots,x_j\rangle$ for some $1\le j\le n$. In fact, let $P\in Ass(\ol{I})$ and by \cite[Corollary 1.3.10]{HH}, there exists a monomial $v\in Mon(S)\setminus P$ such that $P=\ol{I}:v$. By \cite[Theorem 1.4.2]{HH},  $v\not\in \ol{I}$ implies that $v^r\not\in I^r$ holds for all integer $r$. For any $x_m\in P$, clearly
$vx_m\in \ol{I}$ holds, thus there exists a positive integer $k$ such that $x_m^kv^k\in I^k$. Since $I^k$ is also monomial of Borel type and $x_m^k\mid x_m^kv^k\in I^k$ holds, thus for any $1\le j< m$, by \cite[Proposition 4.2.9(2)]{HH}, there exists an integer $t\ge 0$ such that $x_j^tv^k\in I^k$. $t>0$ holds since $v^k\not\in I^k$. Then $x_j^tv^{k-1}\in I^k:v\seq \ol{I}:v=P$. By the choice of $v$, we have $v\not\in P$ thus $v^{k-1}\not\in P$. Then $x_j\in P$ and it shows that $\ol{I}$ is of Borel type.

(2) Now assume that $I$ is strongly stable. Then for any monomial $u\in \ol{I}$, there exists a positive integer $k$  such that
$u^k\in I^k$. If  $x_j\mid u$, then $x_j^k\mid u^k$. Assume $u^k=w_1w_2\cdots w_k$, in which  $w_i\in Mon(S)\cap I$. Assume further that $x_j^{a_i}\mid w_i$, where $\sum_{i=1}^ka_i=k$ and $a_i\ge 0$. Then for any $i<j$, we have $x_i^{a_i}(w_1/x_j^{a_i})\in I$.
Then
$$[x_i(u/x_j)]^k=\prod_{i=1}^k x_i^{a_i}(w_i/x_j^{a_i})\in I^k,$$
thus by \cite[Theorem 1.4.2]{HH}, $x_i(u/x_j)\in \ol{I}$ holds. This shows that $\ol{I}$ is strongly stable.

(3) Assume that $I$ is Borel-fixed. Just as in $(2)$, we assume $u^k = w_1w_2\cdots w_k \in I^k$, where $w_i\in I$. Note that $(\alpha(u))^k = \alpha(u^k)$, it will suffice to show that $\alpha(u^k) \in I^k$ for every $\alpha \in \mathcal{B}$, where $\mathcal{B}$ is the set of upper invertible $n\times n$ matrices over $K$. By Proposition \ref{operation}(2), it is clear since $I$ is Borel-fixed.

(4) Assume that $I$ is lexsegment. For each $u \in \ol{I}$, there exists a positive integer $k$, such that $u^k = \prod_{l=1}^ku_l \in I^k$. Let $u=x_i^{a_i}(\prod_{j=1}^{i-1}x_j^{a_j})(\prod_{t=i+1}^{n}x_t^{a_t})$, and let $v=x_i^{b_i}(\prod_{j=1}^{i-1}x_j^{a_j})(\prod_{t=i+1}^{n}x_t^{b_t})$  such that $b_i > a_i$ and $\sum_{t=i}^n b_{t}=\sum_{t=i}^n a_{t}$. Assume that $u_l=x_i^{a_{li}}(\prod_{j=1}^{i-1}x_j^{a_{lj}})(\prod_{t=i+1}^{n}x_t^{a_{lt}})$ for $1 \le l \le k$. It is easy to see $\sum_{l=1}^k a_{lj}=ka_j$ for $1 \le l \le n$. In the following, we will show that there exist $v_1, \ldots, v_k \in I$ such that $v^k = \prod_{l=1}^k v_l$, which implies $v \in \ol{I}$. In fact, we can choose $v_l$ under the following rule: If $\prod_{t=i+1}^{n}x_t^{a_{lt}}=1$, then set $v_l=u_l$ and $v_l' = 1$; If $\prod_{t=i+1}^{n}x_t^{a_{lt}} \neq 1$, then set $v_l = x_i^{a_{li}+1}(\prod_{j=1}^{i-1}x_j^{a_{lj}}) \cdot v_l'$, such that $degree(v_l') = degree(u_l) - \sum_{j=1}^i a_{lj}-1$ and $v_l' \,|\, \prod_{t=i}^n x_i^{b_i}/\prod_{t=1}^{l-1} v_t'$ with the exponent of $x_i$ as small as possible. Note that $a_i < b_i$ and $degree(v)=degree(u)$, there exist a group of $v_1, \ldots, v_k$ such that $v^k = \prod_{l=1}^k v_l$.

(5) Assume that $I$ is universal lexsegment. If the minimal generating set of $I$ is $G_{min}(I)= \{ u_1, \ldots, u_m \}$ with $u_1 > u_2 > \cdots > u_m$ by pure lexicographic order, then there exists a group of positive integers $a_1, \ldots, a_m$, such that $u_i= x_i^{a_i} \prod_{j=1}^{i-1}x_j^{a_j-1}$ for each $1 \le i \le m$. Let $\mathcal{C}(I)$ be the convex hull of the set of lattice points $\{ \alpha \,|\, x^{\alpha} \in I\}$. By Corollary 1.4.3 \cite{HH}, $\ol{I}= \langle x^{\alpha} \,|\, \alpha \in \mathcal{C}(I) \rangle$. Note that the structure of $u_i$ for $1 \leq i \leq m$, it is not hard to see that $\ol{I} = I$. Hence $\ol{I}$ is universal lexsegment.
\end{prof}

\vs{3mm} We remark that Theorem  \ref{Closure} (1) can also be proved in a similar way as is used in proving $(2)$ and $(3)$.

It is known that $I\seq \ol{I}\seq \sqrt{I}$ holds for every ideal of a (noetherian) ring $R$. Thus $\sqrt{I}=\sqrt{\ol{I}}$ holds. By the primary decomposition theorem (see \cite{AM, Eisenbud}), we record

\begin{prop}\label{Minimal Prime ideals} For any ideal $I$ of a noetherian ring $R$, $Min(I)=Min(\ol{I})$  holds, where $Min(I)$ is the set of all
prime ideals minimal over $I$. In particular, a squarefree monomial ideal $I$ of $S$ is integrally closed.
\end{prop}

\begin{prop}\label{integral power inclusion} For a monomial ideal $I$ of $S$ and any integer $k \geq 1$, $\ol{I}^k \seq \ol{I^k}$ holds.
\end{prop}

\begin{prof} First, note that $$\ol{I^k} = \langle \{ u\in S \mid \exists\, l\,\,such\,\, that\,\, u^{l} \in I^{kl}\} \rangle$$ and
$$\ol{I}^k = \langle \{ \prod_{i=1}^k w_i \,|\, \exists l_i,\,\, such \,\,that\,\, w_i^{l_i} \in I^{l_i} \} \rangle.$$
For every $v = \prod_{i=1}^k w_i \in \ol{I}^k$ with $w_i^{l_i} \in I^{l_i}\,(\forall i=1, \ldots, k)$, let $l=lcm(l_1, \ldots, l_k)$. Then $w_i^{l} \in I^{l}$ holds for each $i=1, \ldots, k$. Thus $v^l = \prod_{i=1}^k w_i^l \in I^{kl}$, which implies $v \in \ol{I^k}$.
\end{prof}

\vs{3mm} The converse inclusion does not hold even for squarefree monomial ideals. We include a counterexample below:
\begin{ex}\label{not equal} Let $u=\prod_{i=1}^6x_i$, and let $$I = \langle x_1x_2x_3, \,x_1x_4x_5,\, x_2x_4x_6,\, x_3x_5x_6\rangle$$ be a squarefree monomial ideal of $S=K[x_1, x_2, x_3, x_4, x_5, x_6]$, thus $I=\ol{I}$. It is easy to check  $u \notin I^2$, but $u^2 \in (I^2)^2$ holds and hence $u \in \ol{I^2}$. Thus $\ol{I^2}\not\seq \ol{I}^2.$
\end{ex}

\section{The $k^{th}$ symbolic power $I^{(k)}$ of an ideal $I$}

Let $I$ be any ideal of a noetherian ring $R$. It follows from \cite[Corllary 2.19]{Eisenbud} that $Min(I)=Min(I^k)$ holds for all positive integer $k$, thus $\underset{k\ge 1} {\cup} Min(I^k)=Min(I)$. Recall that for each $P\in Min(I)$, $ker(R\mapsto (R/I)_P)$ is the $P$-primary component of $I$, and it depends only on $I$ and $P$ in an irredundant primary decomposition
of $I$. If $$I^k=\underset{P\in Ass(I^k)}{\cap}Q(P) $$ is an irredundant primary decomposition of $I^k$, then $Q(P)=ker(R\to (R/I^k)_P)$
holds for each $P\in Min(I^k)$, and $\underset{P\in Min(I^k)}{\cap}Q(P)$ is independent of the primary decomposition of $I^k$. Recall that
$$I^{(k)}=\underset{P\in Min(I)}{\cap}ker(R\to (R/I^k)_P)$$ is called the {\it $k^{th}$ symbolic power} of $I$.

By \cite[Section 3]{HTT}, $$I^{(k)}=I^k:(\cap_{P\in Ass^*(I)\setminus Min(I)} P)^\infty,$$ where $Ass^*(I)=\underset{k\ge 0} {\cup} Ass(I^k)$. Thus it follows from Proposition \ref{operation} that if $I$ is monomial of Borel type, then so is $I^{(k)}$. In the following, we will give a direct and alternative proof to the fact. We need some preparations.

\begin{mdef}\label{new graded} {\it Let  $B$ be a nonempty subset of $[n]$. For any monomial $u= x_1^{a_1}x_2^{a_2}\cdots x_n^{a_n}$, $\underset{j\in B}{\sum}{a_{j}}$ is called {\it the $B$-degree} of the monomial $u$. An ideal $I$ is called $B$-graded if $f_i\in I$ holds for the $B$-graded decomposition  $f=f_0 + f_1 + \cdots + f_k$ of each $f\in I$, where $f_i$ is the $B$-degree $i$ component of $f$.}

\end{mdef}

It is easy to check the following property.

\begin{lem}\label{monomial to graded} If $I$ is a monomial ideal, then for every subset $B$ of $[n]$, $I$ is $B$-graded.
\end{lem}

\vs{3mm}Let $A$ be a subset of $[n]$. For a monomial $u=x^{\al} \in S$ with $\al = (a_1, \ldots, a_n)$, denote $\al(A)= (b_1, \cdots, b_n)$, where
\begin{equation}
\begin{cases}
b_i = a_i\,\quad if\,\, i\in A, \\
b_i = 0\,\quad if\,\, i\in [n]\setminus A.
\end{cases}
\end{equation}
Denote $u(A)= x^{\alpha(A)}$. We also denote $M(A)= \{u(A) \,|\, u \in M\}$ for any nonempty subset $M$ of $Mon(S)$.

For a prime ideal $P$ and an ideal $I$ of $S$,  denote
$$J(I,P) = \{f\in S \,|\, \exists\, g\in S\setminus P,\,such\, \,that\,\,  fg\in I\}.$$
Note that $J(I,P)=ker(S\to (S/I)_P).$
For  a monomial ideal $I$, let $G(I)$ ($G_{min}(I)$) be its (minimal) generating set of monomials, and denote $I(A)=G_{min}(I)(A)$ for a subset $A$ of $[n]$.

\begin{prop}\label{Borel evaluation}
Let $I$ and $P$ be monomial ideals of $S$. If $P$ is a prime ideal, then $J(I,P)$  is a monomial ideal. Furthermore,
$$I(X_P)=\{u(X_P)\,|\,u\in G_{min}(I)\}$$
is a monomial generating set of $J(I,P)$, where $X_P=\{i\in [n] \,|\, x_i \in P\}$. In particular, $$|G_{min}(J(I,P))| \leq |G_{min}(I)|.$$
\end{prop}

\begin{prof} Assume $P=\langle x_{i_1}, \cdots, x_{i_t}\rangle$ and let $X_P=\{i_1,\ldots,i_t\}$. First, we will show that $J(I,P)$ is a monomial ideal.
For any $f\in J(I,P)$ and any $g\in S\setminus P$ such that $fg\in I$, let $f=f_0 + f_1 + \cdots + f_m$ and $g=g_0 + g_1 + \cdots + g_l$ be their $X_P$-graded decompositions. Then
 $$(fg)_r= \sum_{i+j=r}f_ig_j.$$
It follows that $f_0g_0=(fg)_0\in I$ since $I$ is a $X_P$-graded ideal by Lemma \ref{monomial to graded}.

We will prove  that $Supp(f) \subseteq J(I,P)$ holds by induction on the graded component number $m$ of $f$. Since $Supp(f_0)\seq  Supp(f)$, it will suffice to show that $supp(f_0) \subseteq J(I,P)$ holds. For this purpose, let $f_0 = u_1 + \cdots + u_s$ and $g_0 = v_1 + \cdots + v_c$, where $$Supp(f_0)=\{u_i\,|\,1\le i\le s\},\,\, Supp(g_0)=\{v_j\,|\, 1\le j\le c\},$$ $u_1<u_2<\cdots<u_s$ and $v_j<v_{j+1}$ under a suitable monomial order. Then it follows that $u_1v_1 \in I$ holds since $I$ is a monomial ideal.

Note that $v_1$ has degree 0 under the $X_P$-grading, thus $v_1\not\in P$ and hence $u_1v_1 \in I$ implies $u_1\in J(I,P)$. Since
$$ u_1v_1g_0 + (u_2 + \cdots + u_s)v_1g_0=f_0g_0v_1 \in I,$$  it follows that $(u_2 + \cdots + u_s)v_1g_0 \in I$. 
Note that both $v_1$ and $g_0$ have degree 0 under the $X_P$-grading, so does $v_1g_0$.
It follows by induction that $Supp(f_0)\seq J(I,P)$ holds. This proves
that $J(I,P)$ is a monomial ideal.

For the second statement, for a monomial $u\in I(X_P)$, there exists a $v \in G_{min}(I)$, such that $u = v(X_P)$. Note that $v(X_P)v([n]\setminus X_P) = v \in I$, and $v([n]\setminus X_P) \in S\setminus P$, so $u=v(X_P) \in J(I,P)$. On the other hand, if a monomial $u\in J(I,P)$, then there exists a monomial $w\in S\setminus P$, such that $uw\in I$. Note that $u(X_P)=(uw)(X_P)$, there exits a monomial $v\in G_{min}(I)$, such that $v | uw$, and hence $v(X_P) | u(X_P)$. Thus $J(I,P)$ is generated by $I(X_P)$.

The last statement is clear.\end{prof}

\begin{cor}\label{disappear large variety} Let $P = \langle x_{i_1}, \cdots, x_{i_k}\rangle$ with $x_j\notin P$ and $x_t \in P $ for every $t < j$. If a monomial ideal $I$ is of Borel type, then for every monomial $u \in G_{min}(J(I,P))$, $x_l \nmid u$ for each $l \geq j$. In particular, if $x_1 \not\in  P$, then $J(I,P) = S$.
\end{cor}

\begin{prof} Let $B = \{1, \cdots, j-1\}$.
It will suffice to show that for every $u \in I$, there exists a $t \geq 0$, such that $u(B)x_j^t \in I$. But this is easy to check whenever $I$ is  of Borel type.
\end{prof}

\vs{3mm}Note that the above conclusion is still true when $I$ is Borel-fixed, strongly stable,  lexsegment, or universal lexsegment.




\vs{3mm}By Proposition \ref{Borel evaluation} and  Corollary \ref{disappear large variety}, the following proposition can be checked directly, so we omit part of the proof.

\begin{prop}\label{Borel transposition} Let $P$ be a monomial prime ideal. If $I$ is of Borel type (strongly stable,  Borel-fixed, lexsegment, or universal lexsegment respectively), then $J(I,P)$ is  of Borel type (strongly stable,  Borel-fixed, lexsegment, or universal lexsegment respectively).
\end{prop}

\begin{prof} We only prove the case when $I$ is of Borel type.
Let $P = \langle x_{i_1}, \cdots, x_{i_k}\rangle$ with $x_j\notin P$ and $x_t \in P $ for every $t < j$. Denote $X_P=\{t \,|\, x_t \in P\}$ and $B=\{1, \cdots, j-1\}$. Clearly, $B \subseteq X_P$.
For a monomial $u \in G_{min}(J(I,P))$, by Corollary \ref{disappear large variety} and the definition of $J(I, P)$, there exists a monomial $w \in S\setminus P$, such that $uw \in I$ and $(uw)(B)=u$. For every pair of $m<l$, if $x_l | u$, then there exists $a \geq 0$ such that $x_m^a(uw/x_l) \in I$, since $I$ is of Borel type. Let $y=x_m^a(uw/x_l)$ and note that $$y(X_P)=y(B)=x_m^a(u/x_l),$$ hence $x_m^a(u/x_l) \in J(I,P)$.
\end{prof}

\vs{3mm}Note that for a universal lexsegment ideal $I$, $depth(S/I)=n-|G_{min}(I)|$, see \cite{MH}. By Proposition \ref{Borel transposition}, $J(I, P)$ is also universal lexsegment. In order to consider the depth of $S/{J(I,P)}$, we need $J(I, P)$ to be a proper ideal of $S$.

\begin{lem}\label{prime inclusion} For a monomial ideal $I$  and a monomial prime ideal $P$ of $S$,  $I \seq P$  holds if and only if $I(X_P)$ generates a proper ideal of $S$, i.e., $J(I,P)$ is a proper ideal of $S$.
\end{lem}

\begin{prof}  If  $P=\langle x_{i_1}, \cdots, x_{i_k}\rangle$ and  is prime over $I$, then for each monomial $u \in I \subseteq P$, $u(X_P) \neq 1$, hence $\langle I(X_P)\rangle \neq S$. On the other hand, if a prime ideal $Q$ does not contain $I$, then there exists a monomial $v\in I \setminus Q$, such that $x_j \nmid v$ for every $j \in X_Q$. Thus $v(X_Q) = 1$, and hence $\langle I(X_Q)\rangle = S$. This completes the proof.
\end{prof}

\vs{3mm}By Proposition \ref{Borel evaluation}, Proposition \ref{Borel transposition} and Lemma \ref{prime inclusion}, the following corollary is direct to check, so we omit the proof.

\begin{cor}\label{depth} Let $I$ be a monomial ideal, and  $P$ a monomial prime ideal containing $I$. If further $I$ is universal lexsegment, then $depth(S/I) \leq depth(S/{J(I,P)})$ holds. Furthermore, the identity holds true if and only if $\{ x_1, \ldots, x_{|G_{min}(I)|}\} \seq P$.
\end{cor}

\vs{3mm}Now we are ready to prove the afore mentioned result:

\begin{thm}\label{SymbolicPower} Let $I$ be a monomial ideal of $S$. If $I$ is strongly stable (Borel-fixed, or of Borel type, respectively), then the $k^{th}$ symbolic power $I^{(k)}$
is also a monomial ideal which is strongly stable (Borel-fixed, or of Borel type, respectively).
\end{thm}

\begin{prof} First we claim that $I^{(k)}$ is a monomial ideal. This can follow from the primary decomposition theorem (see e.g., \cite[Theorem 3.10]{Eisenbud}), together with \cite[Theorem 1.3.1 and Proposition 1.3.7]{HH}. For any $P\in Min(I)$, note also that
$$ker(R\to (R/I^k)_P )= J(I^k,P),$$ thus  gives a direct proof  to the fact.

Note that $I^k$ is strongly stable (Borel-fixed, or of Borel type, respectively), if $I$ is strongly stable (Borel-fixed, or of Borel type, respectively). Hence for every $P\in Min(I)$, $ker(R\to (R/I^k)_P) = J(I^k,P)$ implies that it is strongly stable (Borel-fixed, or of Borel type, respectively) by Proposition \ref{Borel transposition}. Thus $$I^{(k)}=\cap_{P\in Min(I)}ker(R\to (R/I^k)_P)$$ is strongly stable (Borel-fixed, or of Borel type, respectively) by Proposition \ref{operation}.
\end{prof}



\vs{3mm} We remark that for a lexsegment ideal $I$, $I^{(k)}$ may be not lexsegment.


\begin{ex}\label{u not p} Let $S= K[x_1, x_2, x_3]$, and let $I= \langle x_1, x_2 \rangle$. Clearly, $I$ is universal lexsegment. But $I^{(2)} = I^2= \langle x_1^2, x_1x_2, x_2^2 \rangle$ is not lexsegment.
\end{ex}


\vs{3mm} In the following, we will show that there exists a class of ideals whose symbolic powers are lexsegment.

\begin{mdef}\label{power lexsegment} A monomial ideal $I$ is called stably lexsegment ideal, if $I^k$ is lexsegment for each $k > 0$.
\end{mdef}

\vs{3mm} Example \ref{u not p} also shows that a universal lexsegment ideal may be not stably  lexsegment. The following proposition shows that a stably lexsegment ideal may be not universal lexsegment. We omit the proof.

\begin{prop}\label{two invariants p l} Let $S = K[x_1, x_2]$. If $I$ is a lexsegment ideal of $S$, then $I$ is stably lexsegment.
\end{prop}

Even though $I$ being lexsegment does not imply $I^{(k)}$ being lexsegment, we have the following conclusion.



\begin{prop}\label{symbolic lexsegment} If $I$ is a stably lexsegment ideal of $S$, then $I^{(k)}$ is lexsegment for each positive integer $k$.
\end{prop}

\begin{prof} If $I$ is stably lexsegment, then $I^k$ is lexsegment for each positive integer $k$. By Proposition \ref{Borel transposition}, for each monomial prime ideal $P$, $J(I^k, P)$ is lexsegment. By Proposition \ref{operation}, $I^{(k)}=\cap_{P\in Min(I)}J(I^k, P)$ is lexsegment.
\end{prof}

\vs{3mm} We end this section with a general result on $I^{(k)}$ for an ideal $I$ in a Noetherian ring, and will improve the result for monomial ideals in section 4.

\begin{prop}\label{SymbolicPower2} Let $I$ be any ideal of a noetherian ring $R$. For each $P\in Min(I)$, let $Q(P)$ be the primary component of the isolated prime ideal $P$ of $I$. Then

$(1)$ $\underset{P\in Min(I)}{\cap}Q(P)^k\seq I^{(k)}\seq \underset{P\in Min(I)}{\cap} Q(P)^{(k)}$
holds true.

$(2)$ If $Q(P)^k$ is primary for each $P\in Min(I)$, then $I^{(k)}=\underset{P\in Min(I)}{\cap}Q(P)^k.$

$(3)$ For any positive integer $k$, $I^{(k)}=I^k$ if and only if  $Ass(I^k)\seq Ass(I)$ holds.

\end{prop}

\begin{prof}  (1) Let $I=\underset{P\in Ass(I)}{\cap}Q(P)$ be any irredundant primary decomposition of $I$. For any $P\in Min(I)$ and any $P_1\in Ass(I)\setminus\{P\}$, there exists an element
$u\in P_1\setminus P$. Thus $u^m\in Q(P_1)\setminus P$ holds for some $m\ge 1$. Thus $Q(P_1)S_P=S_P$ and hence $IS_P=Q(P)S_P$. Then $I^kS_P=Q(P)^kS_P$ and hence $Q(P)^k\seq Ker(R\to (R/I^k)_P)$. This shows $$\underset{P\in Min(I)}{\cap}Q(P)^k\seq I^{(k)}.$$

On the other hand, for any $r\in Ker(R\to (R/I^k)_P)$, there exists an element $s\not\in P$ such that $sr\in I^k$. Then $sr\in Q(P)^k \seq Q(P)^{(k)}$. Since $Q(P)^{(k)}$ is $P$-primary, it follows that $r\in Q(P)^{(k)}$. This shows $ I^{(k)}\seq \underset{P\in Min(I)}{\cap} Q(P)^{(k)}.$

(2) This follows from (1).

(3) Consider an irredundant primary decomposition
$$I^k=\underset{P\in Ass(I^k)}{\cap}Q(P)$$
of $I^k.$ Since $Min(I)=Min(I^k)\seq Ass(I^k)$ always holds and
$$I^{(k)}=\underset{P\in Min(I^k)}{\cap}Q(P),$$
it follows that $I^k=I^{(k)}$ holds if and only if $Min(I^k)=Ass(I^k)$, and the latter holds if and only if $Ass(I^k)\seq Ass(I)$.

\end{prof}

\vs{3mm}We remark that $Min(Q(P)^{k})=\{P\}$ and in fact $$Q(P)^{(k)}=ker[R\to (R/Q(P)^k)_P],$$ see \cite[Theorem 23, page 232]{Zariski}.

\begin{cor}\label{Squarefree}  (\cite[Proposition 1.4.4]{HH}) Let $I$ be a squarefree  monomial ideal of a polynomial ring $K[x_1,\ldots,x_n]$ over a field $K$. Then for any $k$,
$$I^{(k)}=\underset{P\in Min(I)}{\cap}P^k.$$

\end{cor}

\begin{prof} For a squarefree  monomials ideal $I$, $Ass(I)=Min(I)$ and for each $P\in Min(I)$, we have $Q(P)=P=\langle x_{i_1},\ldots, x_{i_r}\rangle$. Then it is easy to verify that $\langle x_{i_1},\ldots, x_{i_r}\rangle^k$ is $P$-primary, and the result follows from Proposition \ref{SymbolicPower2}(2). \end{prof}

\section{A simplicial complex and a decomposition of $I^{(k)}$ for a monomial ideal $I$}

In this section, we will use notations established before  to improve Proposition \ref{SymbolicPower2} (1) for monomial ideals. In doing so, we will define and study a new simplicial complex.

\vs{3mm} For a subset $A$ of $[n]$, $\langle I(A)\rangle = S$ if and only if $1 \in G_{min}(I)(A)$, and the latter holds if and only if there exists a monomial $x^{\al}=x_1^{a_1}\cdots x_n^{a_n}\in G_{min}(I)$ in which $a_i\neq 0$ implies $i\notin A$. Thus if $B\subseteq A$ and $\langle I(A)\rangle = S$, then clearly $\langle I(B)\rangle =S$ also holds.

\begin{mdef}\label{eliminating simplicial complex} For any monomial ideal $I$ of $S$, there is the following simplicial complex $$_I\bigtriangleup = \{A\subseteq [n] \,|\, \langle I(A)\rangle = S\}.$$
It will be called {\it the eliminating simplicial complex of $I$}.

\end{mdef}

We remark that a simplicial complex on $[n]$ usually contains all the singletons, but we do not assume this condition.
By Lemma \ref{prime inclusion}, it is easy to prove the following proposition.

\begin{prop}\label{minimal primes} If $I$ is a monomial ideal of $S$, then $Min(I) = \{P_B \,|\, B \in \mathcal{N}(_I\bigtriangleup)\}$, where $\mathcal{N}(_I\bigtriangleup)$ consists of the minimal nonfaces of $_I\bigtriangleup$ and $P_B=\langle \{x_i\,|\,i\in B\} \rangle$.
\end{prop}

In the following, we will consider about the radical ideal $\sqrt{I}$ of $I$.
Note that for a monomial ideal $u$ and a subset $B$ of $[n]$, $u(B) = 1$ if and only if $\sqrt{u}(B)=1$. So Lemma \ref{prime inclusion} implies the following well known property.
\begin{cor}\label{radical min} If $I$ is a monomial ideal of $S$, then $_{\sqrt{I}}\bigtriangleup = _I \bigtriangleup$. In particular, $Min(\sqrt{I}) = Min(I)$.
\end{cor}

Let $J$ be a monomial ideal of $S$. Recall from \cite{Villareal, HH} that the Stanley-Reisner ideal of the simplicial complex $_J\bigtriangleup$ is the ideal $I_{_J\bigtriangleup}$ of $S$, which is generated by the squarefree monomials $x_B= \prod_{i\in B} x_i$ with $B \notin\, _J\bigtriangleup$. The Alexander dual of $_J\bigtriangleup$, denoted by $_J\bigtriangleup^{\vee}$, is defined by $_J\bigtriangleup^{\vee} = \{[n]\setminus B \,|\, B\notin _J\bigtriangleup\}$. It is easy to see that for a subset $B$ of $[n]$, $B \in \mathcal{N}(_J\bigtriangleup)$ if and only if $[n]\setminus B \in \mathcal{F}(_J\bigtriangleup^{\vee})$. We have the following observation.

\begin{prop}\label{Alexander dual} If $J$ is a monomial ideal of $S$, then $I_{_J\bigtriangleup^{\vee}} = \sqrt{J}$. In particular, if $J$ is a squarefree monomial ideal, then $I_{_J\bigtriangleup^{\vee}} = J$.

\end{prop}

\begin{prof} Note that $\sqrt{J}$ is squarefree, so $\sqrt{J}=\cap_{P\in Min(\sqrt{J})} P = \cap_{P\in Min(J)} P$ by Corollary \ref{radical min}. By Proposition \ref{minimal primes}, $Min(J) = \{P_B \,|\, B \in \mathcal{N}(_J\bigtriangleup)\}$.  Hence the standard primary decomposition of $I_{_J\bigtriangleup^{\vee}}$ is
$$I_{_J\bigtriangleup^{\vee}}=\underset{B\in \mathcal{F}(_J\bigtriangleup^{\vee})}{\cap} P_{[n]\setminus B}=\underset{A\in \mathcal{N}(_J\bigtriangleup)}{\cap} P_{A}=\underset{P\in Min(J)}{\cap} P = \sqrt{J}.$$

Note that $\sqrt{J} = J$ while $J$ is squarefree, so the second part is clear.
\end{prof}

\vs{3mm}For a monomial $u=x^{\alpha}$ with $\alpha = (a_1 , \cdots , a_n)$, denote $A(u)=\{i \,|\, a_i \neq 0\}$. In the following, we will show that the inclusions  appeared in Proposition \ref{SymbolicPower2}(1) are actually equalities for a monomial ideal $I$ of $S$. For this purpose, we need the following Lemmas.

\begin{lem}\label{evaluation identity} Let $I$ be a monomial ideal of $S$, and $B$ a subset of $[n]$. If $A(u) \seq B \seq [n]$ holds for each $u \in G_{min}(I)$, then $J(I,P_B) = I$.
\end{lem}

\begin{prof} By Proposition \ref{Borel evaluation}, $J(I,P_B) = \langle I(B)\rangle$ holds. Note also that for each $u \in G_{min}(I)$, $A(u) \seq B$ holds by assumption, hence $u(B) = u$ holds for every $u \in G_{min}(I)$. This is equivalent to saying that  $J(I, P_B) = I$.
\end{prof}

\vs{3mm} In the following lemma, let $G^k=\{\prod_{i=1}^ku_i\,|\,u_i\in G\}$ and set $G^k(B)=(G^k)(B)$.

\begin{lem}\label{power evaluation identity} Let $G$ be a set of monomials of $S$, and Let $B$ be a subset of $[n]$. Then for any positive integer $k$, the identity $\langle G(B)\rangle^k = \langle G^k(B)\rangle$ holds. 
In particular, $\langle G^k(B)\rangle =S$ holds if and only if $\langle G(B)\rangle =S$.
\end{lem}

\begin{prof} Clearly, we only need to prove the first statement. By definition, both $\langle G(B)\rangle^k$ and $\langle G^k(B)\rangle$ are monomial ideals of $S$. On the one hand, for every monomial $u \in \langle G(B)\rangle^k$, there exists $u_1, \cdots, u_k \in G$, such that $\prod_{i=1}^k u_i(B) \,|\, u$. Note that
$$\prod_{i=1}^k u_i(B) = (\prod_{i=1}^k u_i)(B),\quad \prod_{i=1}^k u_i \in G^k,$$
hold, hence $u\in \langle G^k(B)\rangle$. This shows $\langle G(B)\rangle^k \seq \langle G^k(B)\rangle$. The other inclusion follows from a similar argument.
\end{prof}

\vs{3mm} The following result improved Proposition \ref{SymbolicPower2} (1) for a monomial ideal $I$ of $S$. It also follows from \cite[Lemma 3.1]{HTT}. Below we include a direct and detailed proof.

\begin{thm}\label{SymbolicPower3} If $I$ is a monomial ideal of $S$, then for any positive integer $k$, the $k^{th}$ symbolic power is
$$I^{(k)}=\underset{B\in \mathcal{N}(_I\bigtriangleup)}{\cap}J(I^k,P_B).$$
Furthermore, $J(I,P_B)^k= J(I,P_B)^{(k)}$ holds for each $B\in \mathcal{N}(_I\bigtriangleup)$ and thus
$$I^{(k)} = \underset{B\in \mathcal{N}(_I\bigtriangleup)}{\cap}J(I,P_B)^k = \underset{P\in Min(I)}{\cap}[ker(S\to (S/I)_P)]^k$$
holds.
\end{thm}

\begin{prof} The first equality follows from Proposition \ref{minimal primes} and the definition of $I^{(k)}$.

 For the remaining equalities, use Proposition \ref{minimal primes} again to have $$Min(I) = \{P_B \,|\, B \in \mathcal{N}(_I\bigtriangleup)\}.$$ So, the $k^{th}$ symbolic power of $I$ is nothing but
 $$I^{(k)}= \underset{B\in \mathcal{N}(_I\bigtriangleup)}{\cap}J(I^k,P_B).$$
 Note that both $I=\langle G(I)\rangle$ and $I^k=\langle \{G(I)\}^k\rangle$ clearly holds, so that
 $$J(I^k,P_B) = \langle \{G(I)\}^k(B)\rangle = \langle \{G(I)(B)\}^k\rangle = J(I, P_B)^k$$
 also holds by Lemma \ref{power evaluation identity}.
  This shows $$I^{(k)} = \underset{B\in \mathcal{N}(_I\bigtriangleup)}{\cap}J(I,P_B)^k .$$

  For the remaining statement, note that $J(I,P_B)$ is $P_B$-primary with $Min(J(I,P_B))=\{P_B\}.$ Hence
  $$ J(I,P_B)^{(k)} = \underset{Q\in Min(J(I,P_B))}{\cap}J(J(I,P_B)^k,Q) = J(J(I,P_B)^k,P_B)$$
  holds. Finally, note that $A(u) \seq B$ holds for every $u\in G_{min}J(I,P_B)$, it follows that $J(J(I,P_B)^k,P_B)$ $ = J(I,P_B)^k$, and hence  $J(I,P_B)^{(k)} = J(I,P_B)^k$ holds. This completes the proof.
\end{prof}

\vs{3mm} In the end, we include an example to illustrate Theorem \ref{SymbolicPower3}:

\begin{ex}\label{identity example} Let $I = \langle x_1^2x_3^2,\, x_1x_2x_3^2\rangle$. The irredundant primary decomposition of $I$ is $I = \langle x_1^2, x_2\rangle \cap \langle x_1\rangle \cap \langle x_3^2\rangle$. Hence $Ass(I) = \{\langle x_1\rangle, \langle x_1, x_2\rangle, \langle x_3\rangle \}$ and $Min(I) = \{\langle x_1\rangle, \, \langle x_3\rangle \}$. It is easy to check that for each $k\geq 1$, $$I^k = \langle \{x_1^{k+i}x_2^{k-i}x_3^{2k} \,|\, i=0,1, \ldots, k \}\rangle = \langle x_1^k\rangle \cap \langle x_3^{2k}\rangle \cap (\cap_{i=1}^k \langle x_1^{k+i}, x_2^{k-i+1}\rangle ).$$  Hence
$$I^{(k)} = \langle x_1^{k}\rangle \cap \langle x_3^{2k}\rangle = \langle x_1\rangle^k \cap \langle x_3^2\rangle^k = \underset{P\in Min(I)}{\cap}J(I,P)^k.$$ It is also clear that
$$\underset{P\in Min(I)}{\cap}J(I,P)^{(k)} = \langle x_1\rangle^{(k)} \cap \langle x_3^2\rangle^{(k)} = \langle x_1\rangle^k \cap \langle x_3^2\rangle^k = I^{(k)}.$$
\end{ex}

\section{Polarization of universal lexsegment monomial ideal}

Let $I$ be a monomial ideal of $S$, and let $<$ be a monomial order on $S= K[x_1,\ldots, x_n]$, such that $x_n < x_{n-1} < \cdots < x_1$. Let $G_{min}(I) = \{u_1, \ldots, u_m\}$ be the minimal generating set of $I$, where $u_i = \prod_{j=1}^n x_j^{a_{ij}}$ for $i=1, \ldots, m$. Let $a_i = max\{a_{1i}, \ldots, a_{mi}\}$. Recall that the polarization of $I$ is a squarefree monomial ideal $T(I) = \langle v_1, \ldots, v_m\rangle$, where $$v_i = \prod_{j=1}^n \prod_{k=1}^{a_{ij}} x_{jk}$$
for $i=1, \ldots, m$. Let $\prec$ be a monomial order on $$T = K[x_{11}, \ldots, x_{1 a_1}, \ldots, \ldots, x_{n 1}, \ldots, x_{n a_n}]$$ such that $x_{ij} \prec x_{kl}$ if $i > k$ or if $i=k$, $j>l$.

If choose $<$ and $\prec$ to be the same kind of monomial order(e.g., lexicographic order, pure lexicographic order or reverse lexicographic order), which satisfies the above convention, then the polarizing process is order-preserving, i.e., for each pair of $u, v \in I$, $u<v$ if and only if $T(u) \prec T(v)$.

Polarization is a powerful tool for studying quite a few important homological and combinatorial invariants, see, e.g., \cite[Corollary 1.6.3]{HH}. But unfortunately, the property of being  Borel type can not be kept in almost all cases after the process of polarization, as the following example shows:

\begin{ex}\label{polarization} Consider the strongly stable monomial ideal $I=\langle x_1^3,x_1^2x_2,x_1x_2^2\rangle$. After polarization, it becomes
$$J=\langle x_{11}x_{12}x_{13}, \, x_{11}x_{12}x_{21},\, x_{11}x_{21}x_{22}\rangle.$$
Since $J$ is monomial and homogeneous, it is clear that $J$ is not a Borel type monomial ideal under any monomial order.
\end{ex}




\begin{thm}\label{Four equal} Let $I$ be a monomial ideal with $G_{min}(I)= \{u_1, \ldots, u_m\}$, where $u_i = \prod_{j=1}^n x_j^{a_{ij}}$ and $u_1 > u_2 > \cdots > u_m$ by pure lexicographic order. Let $$a_j= max\{a_{1j}, \ldots, a_{mj}\}$$ for $j=1, \ldots, n$. Then the following statement are equivalent.

$(1)$ $I$ is universal lexsegment;

$(2)$ $u_i = x_i^{a_i}\prod_{j=1}^{i-1} x_j^{a_j-1}$ for $i=1, \ldots, m$;

$(3)$ $x_i(u/x_{m(u)}^{b_{m(u)}}) \in I$ holds for each $u= \prod_{j=1}^n x_j^{b_{j}} \in I$ and each $i< m(u)$;

$(4)$ For any monomial $u= \prod_{j=1}^n x_j^{b_{j}} \in I$,  $x_i(u/x_j^{b_{j}}) \in I$ holds for each pair $0\leq i<j \leq n$ with $b_j > 0$.
\end{thm}

\vs{3mm}\begin{prof} $(1) \Leftrightarrow (2)$ is well known (see \cite{MH}), and $(2) \Rightarrow (3)$ is clear.

$(3) \Rightarrow (4)$: For any monomial $u= \prod_{k=1}^n x_k^{b_{k}} \in I$ and each pair $0\leq i<j \leq n$ with $b_j > 0$, $x_{m(u)-1}\prod_{k=1}^{m(u)-1} x_k^{b_{k}}=x_{m(u)-1}(u/x_{m(u)}^{b_{m(u)}}) \in I$ by $(3)$. By induction, $x_{j}\prod_{k=1}^{j} x_k^{b_{k}} \in I$. By $(3)$ again, $x_{i}\prod_{k=1}^{j-1} x_k^{b_{k}} \in I$. Note that $x_{i}\prod_{k=1}^{j-1} x_k^{b_{k}} \,|\, x_i(u/x_j^{b_{j}})$, $x_i(u/x_j^{b_{j}}) \in I$.

$(4) \Rightarrow (2)$: It will suffice to show that for each $1\leq k\leq m$, $u_k = x_k^{a_k}\prod_{j=1}^{k-1} x_j^{a_j-1}$ and $a_{ik} = a_k-1$ for each $k+1 \leq i \leq m$. We will prove it by induction. It is clear that $u_1 = x_1^{a_1}$ under the condition $(4)$, and it is easy to see that $a_i=0$ implies that $a_j=0$ for each $j>i$. Note that $u_1 > u_2 > \cdots > u_m$ by pure lexicographic order. Hence we can assume $u_2 = x_1^{a_{21}}x_2^{a_{22}}$ with $a_{21}<a_1$ and $a_{22}>0$. Note that $x_1^{a_{21}+1} \in I$, so $a_{21}=a_1-1$. In a similar way, we get  $a_{i1} = a_1-1$ for every $i = 2, \ldots, m$.
Now assume that the conclusion holds true for all the integers less than $k$, and we are going to show that $u_{k} = x_{k}^{a_{k}}\prod_{j=1}^{k-1} x_j^{a_j-1}$ and $a_{ik} = a_k-1$ holds for each $k+1 \leq i \leq m$.  By inductive assumption,
$$u_i=(\prod_{j=1}^{k-1} x_j^{a_j-1})(\prod_{j=k}^{m} x_j^{a_{ij}})$$
holds for each $k \leq i \leq m$. If assume to the contrary that $a_{ik} < a_k-1$ holds for some $i$,  then $w= x_k^{a_{ik}+1} (\prod_{j=1}^{k-1} x_j^{a_j-1}) \in I$. By the definition of $a_k$, there exists an integer $t \geq k$ such that $x_k^{a_k} \mid u_t$ holds. Hence $w $ properly divides $u_t$,  contradicting $u_t \in G_{min}(I)$. Hence either $a_{ik} = a_k-1$ or $a_{ik} = a_k$ holds for each $k\leq i\leq m$. Note that $u_k > u_{k+1} > \cdots > u_m$, it is easy to check  $u_k = x_k^{a_k}\prod_{j=1}^{k-1} x_j^{a_j-1}$, and that $a_{ik} = a_k-1$ holds for each $k+1 \leq i \leq m$. This completes the verification.
\end{prof}

\begin{rem}\label{universal} Note that $(3)$, as an equivalent description of universal lexsegment ideal, explores the difference between universal lexsegment ideal, strongly stable ideal and the monomial ideal of Borel type. Actually, a universal lexsegment monomial ideal is a kind of {\it super-stable} monomial ideal.
\end{rem}

If $I$ is a universal lexsegment ideal of $S= K[x_1,\ldots, x_n]$ with the minimal generating set $G_{min}(I)$, it is clear that $|G_{min}(I)| \leq  n$. We call a universal lexsegment ideal to be {\it full}, if $|G_{min}(I)| = n$.



In the following, we will consider about the polarization of some class of monomial ideals with respect to $<$, and characterize the monomial ideals which become squarefree strongly stable with respect to $\prec$ after polarization. Recall that a squarefree monomial ideal $I$ is called {\it squarefree strongly stable}, if for each squarefree monomial $u\in I$ and each pair $j<i$ such that
$x_i\mid u$ but $x_j\nmid u$, one has $x_j(u/x_i)\in I$ (see \cite{AHH} or \cite[page 124]{HH}).

\begin{thm}\label{super-stable form} Let $I$ be a monomial ideal with $G_{min}(I)= \{u_1, \ldots, u_m\}$, where $u_i = \prod_{j=1}^n x_j^{a_{ij}}$ and $u_1 > u_2 > \cdots > u_m$ by pure lexicographic order. Let $$a_j= max\{a_{1j}, \ldots, a_{mj}\}$$ for $j=1, \ldots, n$. Then $T(I)$ is squarefree strongly stable,  if $I$ is  universal lexsegment. Further more, if $a_j \neq 1$ holds for each $j=1, \ldots, n$, and for each $1\leq j < m$, there exists an integer $i$ such that $0 < a_{ij} < a_j$, then $T(I)$ is squarefree strongly stable if and only if $I$ is universal lexsegment.
\end{thm}



\vs{3mm}\begin{prof} If $I$ is universal lexsegment, then $u_k = x_k^{a_k}\prod_{j=1}^{k-1} x_j^{a_j-1}$ for $k = 1, \ldots, m$. Hence $$T(u_k) = (\prod_{i=1}^{k-1} \prod_{j=1}^{a_j-1} x_{ij})(\prod_{j=1}^{a_k} x_{kj}).$$ In the case, it is direct to check that $T(I)$ is squarefree strongly stable. This complete the proof of the first statement.

For the second statement, we only need to prove the necessity part.
 It is easy to see that $u_1 = x_1^{a_1}$. We claim that $a_{i1} =a_1-1$ for every $i = 2, \ldots, m$. In fact, if there exists some $u_t$ such that $a_{t1}<a_1-1$, then consider
 $$T(u_t) = (\prod_{j=1}^{a_{t1}} x_{1j})(\prod_{i=2}^{n} \prod_{j=1}^{a_{ti}} x_{ij}) \in T(I).$$
 Since $T(I)$ is squarefree strongly stable, it follows that  $v = x_{1,a_{t1}+1}(T(u_t)/{x_{21}}) \in T(I)$ holds. Note that for each generating element $u$ of $T(I)$, if $x_{ij} | u$, then $x_{i,j-1} | u$, \ldots, $x_{i,1} | u$ hold. Hence $v \in T(I)$ implies $v_1=x_{1,a_{t1}+1}(T(u_t)/{\prod_{j=1}^{a_{t2}} x_{2j}}) \in T(I)$. Note that for any $j$, $x_{2j} \nmid v_1$, thus $x_{2,a_2}(v_1/{x_{31}}) \in T(I)$ whenever $x_{31} \,|\, v_1$. Repeat the discussion above, it follows that $x_{2,a_2}(v_1/{\prod_{j=1}^{a_{t3}} x_{3j}}) \in T(I)$, and hence $v_1/{\prod_{j=1}^{a_{t3}} x_{3j}} \in T(I)$ holds since $a_2 > 1$. By induction, we have $\prod_{j=1}^{a_{t1}+1} x_{1j} \in T(I)$, contradicting $T(u_1) \in G_{min}(T(I))$. Hence $a_{t1} = a_1 - 1$ holds for each $t=2, \ldots, m$. Finally, repeat a discussion used in proving $(4) \Rightarrow (2)$ of theorem \ref{Four equal}, the result follows by induction.
\end{prof}



\vs{3mm}Now assume that $a_j \neq 1$ holds for each $j=1, \ldots, n$. If there exists some $j$ with $a_j > 1$, such that either $a_{ij}=0$ or $a_{ij}=a_{j}$ for each $i=1, \ldots, m$, then the situation will be a little more complicated than the case in Theorem \ref{super-stable form}. In this case, let $W_I = \{j \,|\, a_j \neq 0\}$, and let $A_I=\{j \in W_I \,|\,$ either $\,a_{ij}=a_j \,$ or $\, a_{ij}=0 \,$ for each $i=1, \ldots, m\}$, and let $B_I = W_I \setminus A_I$. Note that if $T(I)$ is squarefree strongly stable, then $W_I= [r]$ holds for some $r\leq n$. In this case, we can decompose $A_I$ into several mutually disjoint subsets  consisting of
successive integers, and denote $A_I = \cup_{t=1}^k A_t$ where
$$A_t = \{j \in Z_{+} \,|\, m_{t-1}+1 \leq j\leq l_t\}.$$
Similarly, let $B_I = \cup_{t=1}^k B_t$, where $B_t = \{j \in Z_{+} \,|\, l_{t}+1 \leq j\leq m_t\}.$

For each $j\in A_I$, either $a_{ij}=0$ or $a_{ij}=a_j$ holds for any $i=1, \ldots, m$. Consider the following two cases:
If $a_{ij}=0$, we claim that it implies $a_{i,j+1}=a_{i,j+2}=\cdots = a_{in}=0$. In fact, if $a_{it} > 0$ for some $t>j$, then $T(u_i)\in T(I)$ holds with $x_{t1} \mid T(u_i)$. Since $T(I)$ is squarefree strongly stable, $x_{j,a_j}(T(u_i)/x_{t1}) \in T(I)$ holds. Note also that $a_j >0$ holds, thus it implies $T(u_i)/x_{t1} \in T(I)$, contradicting $u_i \in G_{min}(I)$. On the other hand, if $a_{ij}=a_j$ holds, then the polarization of $u_i$ contains $\prod_{t=1}^{a_j} x_{jt}$ as its factor, which contains all the indeterminants related to
$x_j$. Note that in each case, essentially there is little change in the problem we are working with.
Due to this reason, the following proposition is routine to check, and we omit the detailed proof.

\begin{prop}\label{two mixed}  If $a_j \neq 1$ for each $j=1, \ldots, n$, then $T(I)$ is squarefree strongly stable if and only if $I = \sum_{i=1}^{k} ( L_iM_i \prod_{j=1}^{i-1} (L_j M_j^{'})),$ where $$L_j= \langle \prod_{t=m_{j-1}+1}^{l_{j}} x_t^{a_t} \rangle,\,\, M_j^{'}=  \langle \prod_{t=l_{j}+1}^{m_{j}} x_t^{a_t-1} \rangle,$$ for every $j=1, \ldots, k$, and $M_j$ is a monomial ideal of  $S$ generated by a full universal lexsegment ideal in the polynomial ring $K[x_{l_j+1}, x_{l_j+2}, \ldots, x_{m_j}]$, where $$0=m_0 < l_1 \leq m_1 < l_2 \leq m_2 < \cdots < l_k \leq m_k \leq n.$$
\end{prop}

\begin{rem}\label{explain} In Proposition \ref{two mixed}, the equality
$$I = \sum_{i=1}^{k} ( L_iM_i \prod_{j=1}^{i-1} (L_j M_j^{'}))$$
can be interpreted in the following: A universal lexsegment ideal on $B_I$ is cut into several
parts by some principal ideals with respect to $A_I$. Note also that the equivalence description would be rather complicated without the assumption $a_j \neq 1$ for each $j=1, \ldots, n$, because it will be some universal lexsegment ideals, squarefree strongly stable ideals being cut into several parts by some principal ideals.
\end{rem}

\begin{ex}\label{universal lexsegment} Let $$I= \langle x_1^3, x_1^2x_2x_3, x_1^2x_2x_4, x_1^2x_2x_5^3x_6^3, x_1^2x_2x_5^3x_6^2x_7^2, x_1^2x_2x_5^3x_6^2x_7x_8^2 \rangle.$$ It is easy to see that $T(I)$ is squarefree strongly stable. By computation, $$a_1=3, a_2=1, a_3=1, a_4=1, a_5=3, a_6=3, a_7=2, a_8=2.$$ An explanation is that $T(I)$ is going to be divided approximately into several parts: for $x_1$, it is universal lexsegment; for $x_2,x_3,x_4$, it is squarefree strongly stable; for $x_5$, it is a principal ideal; for $x_6, x_7, x_8$, again it is universal lexsegment.
\end{ex}






\begin{mdef}\label{exponent vector} Let $I$ be a monomial ideal, and let $G_{min}(I) = \{u_1, \ldots, u_m\}$ with $u_i = \prod_{j=1}^n x_j^{a_{ij}}$ for $i=1, \ldots, m$. The vector $(a_1, a_2, \ldots, a_n)$ is called the {\it exponent vector} of $I$, where $a_i = max\{a_{1i}, \ldots, a_{mi}\}$. For a squarefree monomial ideal $$J \subseteq K[x_{11}, \ldots, x_{1 \tau_1}, \ldots, \ldots, x_{n 1}, \ldots, x_{n \tau_n}],$$ the vector $(b_1, b_2, \ldots, b_n)$ is called the {\it extension  vector} of $J$, where $$b_i = max\{j \,|\, \, \exists \,u \in G_{min}(J)\,\, such \,\,that \,\,x_{ij} \mid u\}.$$
\end{mdef}

By definition, the following Proposition is clear.

\begin{prop}\label{exponent extension equal} For any monomial ideal $I$, The exponent vector of $I$ is equal to the extension vector of $T(I)$.

\end{prop}


\begin{prop}\label{Borel ideals number} Let $J$ be a polarization of some monomial ideal of $S$, with the extension vector $(b_1, \ldots, b_n)$. If $J$ is squarefree strongly stable and $b_i \neq 1$ for each $i=1, \ldots, n$, then the number of monomial $I$ such that $T(I) = J$ is $2^t$, where $t = |W_I|$.

\end{prop}

\begin{prof} By Proposition \ref{two mixed}, $I$ is uniquely determined by the exponent vector $(a_1, \ldots, a_n)$ and the set $A_I$. By Proposition \ref{exponent extension equal}, $(a_1, \ldots, a_n) = (b_1, \ldots, b_n)$ is a constant vector. Hence $I$ is completely determined by the set $A_I$. Note that $A_I$ could be any subset of $W_I$, so the number of $I$ such that $T(I) = J$ is $2^t$, where $t = |W_I|$.
\end{prof}

\end{document}